\documentclass{amsart}
\usepackage{amssymb}

\newtheorem{theorem}{Theorem}

\newtheorem{lemma}{Lemma}

\newcommand{\bee}[1]{\begin{equation}\label{#1}}
\newcommand{\beq}[1]{\begin{eqnarray}\label{#1}}
\newcommand{\ene}{\end{equation}}
\newcommand{\eqe}{\end{eqnarray}}

\begin{document}
\title{On the codimension growth of almost nilpotent Lie algebras}

\author[D. Repov\v s and M. Zaicev]
{Du\v san Repov\v s and Mikhail Zaicev}

%\author{Du\v{s}an~Repov\v{s}}
\address{Du\v san Repov\v s \\Faculty of Mathematics and Physics, and
Faculty of Education, University of Ljubljana,
P.~O.~B. 2964, Ljubljana, 1001 Slovenia}
\email{dusan.repovs@guest.arnes.si}

\address{Mikhail Zaicev \\Department of Algebra\\ Faculty of Mathematics and
Mechanics\\  Moscow State University \\ Moscow, 119992 Russia}
\email{zaicevmv@mail.ru}

\thanks{The first author was partially supported by the Slovenian Research Agency 
grants P1-0292-0101 and J1-2057-0101. 
The second author was partially supported by RFBR  grant No 09-01-00303a}



\begin{abstract}
We study codimension growth of infinite dimensional Lie algebras over a field of
characteristic zero. We prove that if a Lie algebra $L$ is an extension of a nilpotent 
algebra by a finite dimensional semisimple algebra then the PI-exponent of $L$ exists and is a positive integer.
\end{abstract}

\date{\today}

\maketitle

\section{Introduction}\label{intr}

We consider algebras over a field $F$ of characteristic zero.
Given an algebra $A$, we can associate to it the sequence of its codimensions
$\{c_n(A)\}, n=1,2,\ldots~$. If $A$ is an associative algebra with a non-trivial
polynomial identity then $c_n(A)$ is exponentially bounded \cite{Reg72} while
$c_n(A)=n!$ if $A$ is not PI.

For a Lie algebra $L$ the sequence $\{c_n(L)\}$ is in general not exponentially bounded  (see for example, \cite{P}). Nevertheless, a class of Lie algebras
with exponentially bounded codimensions is quite wide. It includes, in particular,
all finite dimensional algebras \cite{B-Dr,GZ-TAMS2010}, Kac-Moody algebras
\cite{Z1,Z2}, infinite dimensional simple Lie algebras of Cartan type
\cite{M1}, Virasoro algebra and many others.

When $\{c_n(A)\}$ is exponentially bounded the upper and the lower limits of
the sequence $\sqrt[n]{c_n(A)}$ exist and the natural question arises: does the ordinary limit $\lim_{n\to\infty} \sqrt[n]{c_n(A)}$ exist? In 1980's Amitsur 
conjectured that for any associative PI algebra such a limit exists and is a 
non-negative integer. This conjecture was confirmed in \cite{GZ1,GZ2}. For Lie 
algebras the series of positive results was obtained for finite dimensional algebras 
\cite{GRZ1,GRZ2,Z3}, for algebras with nilpotent commutator 
subalgebras \cite{PM} and some other classes (see \cite{M2}).

On the other hand it was shown in \cite{ZM} that there exists a Lie algebra $L$ with
$$
3.1 < \liminf_{n\to\infty} \sqrt[n]{c_n(L)} \le 
\limsup_{n\to\infty} \sqrt[n]{c_n(AL)} < 3.9 \quad .
$$
This algebra $L$ is soluble and almost nilpotent, i.e. it contains a nilpotent ideal
of finite codimension. Almost nilpotent Lie algebras are close in some sense
to finite dimensional algebras. For instance, they have the Levi decomposition
under some natural restrictions (see \cite[Theorem 6.4.8]{Baht}), satisfy Capelli
identity, have exponentially bounded codimension growth, etc. Almost nilpotent Lie algebras play an important role in the theory of codimension growth since all minimal
soluble varieties of a finite basic rank with almost polynomial growth are generated by
almost nilpotent Lie algebras. Two of them have exponential growth with ratio $2$
and one is of exponential growth with ratio $3$.

In the present paper we prove the following results.
\medskip

{\bf Teorem 1.}
{\em Let $L$ be an almost nilpotent Lie algebra over a field $F$ of characteristic zero. If
$N$ is the maximal nilpotent ideal of $L$ and $L/N$ is semisimple then the PI-exponent
of $L$,
$$
\mbox{exp}(L)=\lim_{n\to \infty}\sqrt[n]{c_n(L)},
$$
exists and is a positive integer.}

\medskip

Recall that a Lie algebra $L$ is said to be special (or SPI) if it is a Lie subalgebra
of some associative PI-algebra.
\medskip

{\bf Teorem 2.}
{\em Let $L$ be an almost nilpotent soluble special Lie algebra over a field $F$ of characteristic zero. Then the PI-exponent
of $L$ exists and is a positive integer.}
\medskip

Note that the special condition in Theorem \ref{t2} is necessary since the
counterexample constructed in \cite{ZM} is a finitely generated almost nilpotent
soluble Lie algebra satisfying Capelli identity of low rank. Nevertheless, its
PI-exponent $\mbox{exp}(L)$ exists (see \cite{VMZ}) but is not an integer since
$\mbox{exp}(L)\approx 3,6$. Note also that Theorem \ref{t1} generalizes the main result of
\cite{GRZ2} and gives an alternative and easier proof of the integrality of PI-exponent
in the finite dimensional case considered in \cite{GRZ2}.

\section{Preliminaries}

Let $L$ be a Lie algebra over $F$. We shall omit Lie brackets in the product of 
elements of $L$ and write $ab$ instead of $[a,b]$. We shall also denote the right-normed product $a(b(c\cdots d)\ldots)$ as $abc\cdots d$. One can find all basic notions of the theory of identities of Lie algebras  in \cite{Baht}.

Let $\bar F$ be an extension of $F$ and $\bar L=L\otimes_F\bar F$. It is not difficult 
to check that $c_n(\bar L)$ over $\bar F$ coincides with $c_n(L)$ over $F$.
Hence it is sufficient to prove our results only for algebras over an algebraically closed field.

Let now $X$ be a countable set of indeterminates and let $Lie(X)$ be a free Lie 
algebra generated by $X$. Lie polynomial $f=f(x_1,\ldots, x_n)\in Lie(X)$ is an 
identity of a Lie algebra $L$ if $f(a_1,\ldots, a_n)=0$ for any $a_1,\ldots, a_n\in 
L$. It is known that the set of all identities of $L$ forms a T-ideal $Id(L)$ of 
$Lie(X)$, i.e. an ideal stable under all endomorphisms of $Lie(X)$. Denote by 
$P_n=P_n(x_1\ldots,x_n)$ the subspace of all 
multilinear polynomials in $x_1,\ldots,x_n$ in $Lie(X)$. Then the intersection 
$P_n\cap Id(L)$ is the set of all multilinear in $x_1,\ldots,x_n$ identities of $L$.
Since $char~F=0$, the union $(P_1\cap Id(L))\cup (P_2\cap Id(L))\cup\ldots$ completely defines all identities of $L$.

An important numerical invariant of the set of all identities of $L$ is the sequence of codimensions
$$
c_n(L)=\dim P_n(L) \quad {\rm where}\quad P_n(L)=\frac{P_n}{P_n\cap Id(L)},
n=1,2,\ldots~.
$$
If $\{c_n(L)\}$ is exponentially bounded one can define the lower and the upper 
PI-exponents of $L$ as 
$$
\underline{\mbox{exp}}(L) = \liminf_{n\to \infty}\sqrt[n]{c_n(L)}, \
\overline{\mbox{exp}}(L) = \limsup_{n\to \infty}\sqrt[n]{c_n(L)}, 
$$
and
$$
\mbox{exp}(L)= \overline{\mbox{exp}}(L)
=\underline{\mbox{exp}}(L),
$$
the (ordinary) PI-exponent of $L$, in case equality holds.

One of the main tools for studying asymptotics of $\{c_n(L)\}$ is the theory of representations of symmetric group $S_n$ (see \cite{JK} for details). Given a multilinear
polynomial $f=f(x_1,\ldots, x_n)\in P_n$, one can define
\begin{equation}\label{e0}
\sigma f(x_1,\ldots, x_n)=f(x_{\sigma(1)},\ldots, x_{\sigma(n)}).
\end{equation}
Clearly, (\ref{e0}) induces $S_n$-action on $P_n$. Hence $P_n$ is an $FS_n$-module 
and $P_n\cap Id(L)$ is its submodule. Then $P_n(L)=\frac{P_n}{P_n\cap Id(L)}$ is also
an $FS_n$-module. Since $F$ is of characteristic zero, $P_n(L)$ is completely
reducible,
\begin{equation}\label{e1}
P_n(L)=M_1\oplus\cdots\oplus M_t
\end{equation}
where  $M_1,\ldots,M_t$ are irreducible $FS_n$-modules and the number $t$ of summands on the right hand side of (\ref{e1}) is called the $n$th colength of $L$,
$$
l_n(L)=t.
$$
Recall that any irreducible $FS_n$-module is isomorphic to some minimal left ideal of
group algebra $FS_n$ which can be constructed as follows.

Let $\lambda\vdash n$ be a partition of $n$, i.e. $\lambda=(\lambda_1,\ldots,\lambda_k)$
where $\lambda_1\ge\ldots\ge\lambda_k$ are positive integers and
$\lambda_1+\cdots+\lambda_k=n$. The Young diagram $D_\lambda$ corresponding to $\lambda$ is a tableau
$$
D_\lambda=\; \begin{array}{|c|c|c|c|c|c|c|} \hline  &   &
\cdots &  &  & \cdots &   \\ \hline & & \cdots &  \\
\cline{1-4}\vdots  \\ \cline{1-1} \\ \cline{1-1} \end{array},
$$
containing $\lambda_1$ boxes in the first row, $\lambda_2$ boxes in the second row, 
and so on. Young tableau $T_\lambda$ is the Young diagram $D_\lambda$ with the integers $1,2,\ldots,n$ in the boxes. Given Young tableau, denote by $R_{T_\lambda}$ the row stabilizer of $T_\lambda$, i.e. the subgroup of all permutations $\sigma\in S_n$
permuting symbols only inside their rows. Similarly, $C_{T_\lambda}$ is the column stabilizer of $T_\lambda$. Denote
$$
R(T_\lambda)=\sum_{\sigma\in R_{T_\lambda}} \sigma~,\quad
C(T_\lambda)=\sum_{\tau\in C_{T_\lambda}} ({\rm sgn}~\tau)\tau~,\quad
e_{T_\lambda}=R(T_\lambda)C(T_\lambda).
$$
Then $e_{T_\lambda}$ is an essential idempotent of the group algebra $FS_n$ that is
$e_{T_\lambda}^2=\alpha e_{T_\lambda}$ where $\alpha\in F$ is a non-zero scalar. It is
known that $FS_n e_{T_\lambda}$ is an irreducible left $FS_n$-module. We denote its 
character by $\chi_\lambda$. Moreover, if $M$ is an $FS_n$-module with the character
\begin{equation}\label{e2}
\chi(M)=\sum_{\mu\vdash n} m_\mu \chi_\mu
\end{equation}
then $m_\lambda\ne 0$ in (\ref{e2}) for given $\lambda\vdash n$ if and only if
$e_{T_\lambda}M\ne 0$.

If $M=P_n(L)$  for Lie algebra $L$ then the $n$th cocharacter of $L$ is
\begin{equation}\label{e3}
\chi_n(L)=\chi(P_n(L))=\sum_{\lambda\vdash n} m_\lambda \chi_\lambda
\end{equation}
and then
\begin{equation}\label{e4}
l_n(L)=\sum_{\lambda\vdash n} m_\lambda,\quad 
c_n(L)=\sum_{\lambda\vdash n} m_\lambda d_\lambda
\end{equation}
where $m_\lambda$ are as in (\ref{e3}) and
$$
d_\lambda=\deg \chi_\lambda=\dim FS_n e_{T_\lambda}.
$$

Recall that Lie algebra $L$ satisfies Capelli identity of rank $t$ if every 
multilinear polynomial $f(x_1,\ldots,x_n), n\ge t,$ alternating on some 
$x_{i_1},\ldots,x_{i_t}$, $\{i_1\ldots,i_t\}\subseteq\{1,\ldots,n\}$ is an
identity of $L$. It is known (see for example, \cite[Theorem 4.6.1]{GZbook}
that $L$ satisfies Capelli identity of rank $t+1$ if and only if all $m_\lambda$
in (\ref{e4}) are zero as soon as $D_\lambda$ has more than $t$ rows, i.e.
$\lambda_{t+1}\ne 0$.

A useful reduction in the proof of the existence of PI-exponent is given by the
following remark.

\begin{lemma}\label{l1}
Let $L$ be an almost nilpotent Lie algebra with the maximal nilpotent ideal $N$.
Let $\dim L/N=p$ and let $N^q=0$. Then
\begin{itemize}
\item[(1)]
$L$ satisfies Capelli identity of the rank $p+q$; and
\item[(2)]
the colength $l_n(L)$ is a polynomially bounded function of $n$.
\end{itemize}
\end{lemma}

{\em Proof.} Choose an arbitrary basis $e_1,\ldots, e_p$ of $L$ modulo $N$ and an 
arbitrary basis $\{b_\alpha\}$ of $N$. If $f=f(x_1,\ldots,x_n)$ is a multilinear 
polynomial then $f$ is an identity of $L$ if and only if $f$ vanishes under all
evaluations $\{x_1,\ldots,x_n\}\rightarrow B=\{e_1,\ldots,e_p\}\cup \{b_\alpha\}$.

Suppose $n\ge p+q$ and $f$ is alternating on $x_1,\ldots,x_{p+q}$. If 
$\varphi: \{x_1,\ldots,x_n\}\rightarrow B$ is an evaluation such that 
$\varphi(x_i)=\varphi(x_j)$ for some $1\le i<j\le p+q$ then $\varphi(f)=0$ since
$f$ is alternating on $x_i,x_j$. On the other hand if any $e_i$ appears among 
$y_1=\varphi(x_1),\ldots, y_{p+q}=\varphi(x_{p+q})$ at most once then
$\{y_1,\ldots,y_{p+q}\}$ contains at least $q$ basis elements from $N$. Hence
$\varphi(f)=0$ since $N^q=0$ and we have proved the first claim of the lemma. The second assertion now follows from the results of \cite{ZM2}.

\hfill $\Box$

As a consequence of Lemma \ref{l1} we get the following:
\begin{lemma}\label{l2}
If $L$ is an almost nilpotent Lie algebra then the sequence $\{c_n(L)\}$ is
exponentially bounded.
\end{lemma}

{\em Proof.} By Lemma \ref{l1}, there exist an integer $t$ and a polynomial $f(n)$
such that $m_\lambda=0$ in (\ref{e4}) for all $\lambda\vdash n$ with $\lambda_{t+1}
\ne 0$ and $l_n(L)=\sum_{\lambda\vdash n}m_\lambda \le f(n)$. It is well-known
(see for example, \cite[Corollary 4.4.7]{GZbook}) that $d_\lambda=\deg\chi_\lambda
\le t^n$ if $\lambda=(\lambda_1,\ldots,\lambda_k)$ and $k\le n$. Hence we get from (\ref{e4}) an upper bound
$$
c_n(L)\le f(n) t^n
$$
and the proof is completed.

\hfill $\Box$

\section{The upper bound for PI-exponent}

The exponential upper bound for codimensions obtained in Lemma \ref{l2} is not precise. In order to prove the existence and integrality of $\mbox{exp}(L)$ we 
shall find a positive integer $d$ such that $\overline{\mbox{exp}}(L)\le d$ and 
$\underline{\mbox{exp}}(L)\ge d$.

Let $L$ be a Lie algebra with a maximal nilpotent ideal $N$ and a finite dimensional 
semisimple factor-algebra $G=L/N$. Fix a decomposition of $G$ into the sum of simple components
$$
G=G_1\oplus\cdots\oplus G_m
$$
and denote by $\varphi_1,\ldots,\varphi_m$ the canonical projections of $L$ to
$G_1,\ldots,G_m$, respectively. Let $g_1,\ldots,g_k$ be elements of $L$ such that for some $\{i_1,\ldots,i_k\}\subseteq\{1,\ldots,m\}$ one has
$$
\varphi_{i_t}(g_t)\ne 0,\quad \varphi_{i_j}(g_t)=0\quad {\rm for ~all}\quad
j\ne t,~1\le t\le k.
$$
For any non-zero product $M$ of $g_1,\ldots,g_k$ and some $u_1,\ldots,u_t\in N$ we define
the height of $M$ as
$$
ht(M)=\dim G_{i_1}+\cdots+\dim G_{i_k}.
$$
Now we are ready to define a candidate to PI-exponent of $L$ as
\begin{equation}\label{e5}
d=d(L)=\max\{ht(M)\vert 0\ne M\in L\}.
\end{equation}
In order to get an upper bound for $\overline{\mbox{exp}}(L)$ we define the following
multialternating polynomials. Let $Q_{r,k}$ be the set of all polynomials $f$ such that

\begin{itemize}
\item[(1)] $f$ is multilinear, $n=\deg f\ge rk$,
$$
f=f(x_1^1,\ldots,x_r^1,\ldots,x_1^k,\ldots,x_r^k,y_1,\ldots, y_s)
$$
where $rk+s=n$; and
\item[(2)] $f$ is alternating on each set $x_1^i,\ldots,x_r^i, 1\le i\le k$.
\end{itemize}

We shall use the following lemma (see \cite[Lemma 6]{Z3}).

\begin{lemma}\label{l3}
If $f\equiv 0$ is an identity of $L$ for any $f\in Q_{d+1,k}$ for some $d,k$ then
$\overline{\mbox{exp}}(L)\le d$.
\end{lemma}

\hfill $\Box$

Note that Lemma 6 in \cite{Z3} was proved for a finite dimensional Lie algebra $L$.
In fact, it sufficient to assume that $L$ satisfies Capelli identity and that $l_n(L)$ 
is polynomially bounded.

\begin{lemma}\label{l4}
Let $L$ be an almost nilpotent Lie algebra and $d=d(L)$ as defined in (\ref{e5}).
Then $\overline{\mbox{exp}}(L)\le d$.
\end{lemma}

{\em Proof}. Let $N$ be the maximal nilpotent ideal of $L$ and let $N^p=0$. We shall
show that any  polynomial from $Q_{d+1,p}$ is an identity of $L$ and apply Lemma \ref{l3}.

Given $1\le i\le m$, we fix a basis $B_i$ of $L$ modulo $\widetilde G_1+\cdots+
\widetilde G_{i-1}+\widetilde G_{i+1}+\cdots+\widetilde G_m+N$ where $\widetilde G_j$ is 
the full preimage of $G_j$ under the canonical homomorphism $L\to L/N$. In other words, 
$|B_i|=\dim G_i$, $\varphi_i(B_i)$ is a basis of $G_i$ and $\varphi_j(B_i)=0$ for any
$1\le j\ne i\le m$. Fix also a basis $C$ of $N$ and let $B=C\cup B_1\cup\cdots\cup B_m$. Then a multilinear polynomial $f$ is an identity of $L$ if and only if $\varphi(f)=0$ for any evaluation $\varphi: X\to B$.

Suppose $f=f(x_1^1,\ldots,x_{d+1}^1,\ldots, x_1^p,\ldots,x_{d+1}^p, y_1,\ldots,y_s) 
\in Q_{d+1,p}$ is multilinear and alternating on each set $\{x_1^i,\ldots,x_{d+1}^i\}$,
$1\le i \le p$. Given $1\le i\le p$, first consider an evaluation $\rho:X\to L$ such 
that $\rho(x_j^i)= b_j\in B_{t_j}$, $1\le j\le d+1$, with
$$
\dim G_{t_1}+\cdots+ \dim G_{t_{d+1}}\ge d+1
$$
in $L/N=G$. Then by definition of $d$ any product of elements of $B$ containing 
factors 
$b_1,\ldots, b_{d+1}$ is zero, hence $\rho(f)=0$. If $\dim G_{t_1}+\cdots + \dim
G_{t_{d+1}}\le d$ then $b_1,\ldots, b_{d+1}$ are linearly dependent modulo $N$,
say, $b_{d+1}=\alpha_1 b_1+\cdots+\alpha_d b_d+w$, $w\in N$. Then the value of 
$\rho(f)$ is the same as of $\rho':X\to L$ where $\rho'(x_1^i)=\rho(x_1^i)=b_1,
\ldots,\rho'(x_d^i)=\rho(x_d^i)=b_d$, $\rho'(x_{d+1}^i)=w$ since $f$ is alternating
on $x_1^i,\ldots,x_{d+1}^i$. It follows that for any evaluation $\rho: X\to B$ one should take at least one value $\rho(x_j^i)$ in $N$ for any $i=1,\ldots,p$ otherwise 
$\rho(f)=0$. But in this case $\rho(f)$ is also zero since $\rho(f)\in N^p=0$ and we 
have completed the proof.

\hfill $\Box$

\section{The lower bound for PI-exponent}

As in the previous Section let $L$ be an almost nilpotent Lie algebra with the maximal
nilpotent ideal $N$ and suppose that semisimple finite dimensional factor-algebra
$G=L/N=G_1\oplus\cdots\oplus G_m$ where $G_1,\ldots,G_m$ are simple.

\begin{lemma}\label{l5}
Given an algebra $L$ as above, there exist positive integers $q$ and $s$ such that
for any $r=tq, t=1,2,\ldots$, and for any integer $j\ge s$ one can find a multilinear 
polynomial $h_t=h_t(x_1^1,\ldots,x_d^1,\ldots,x_1^r,\ldots,x_d^r,y_1,\ldots, y_{j})$ 
such that:
\begin{itemize}
\item[(1)]
$h_t$ is alternating on each set $\{x_1^i,\ldots,x_d^i\}$, $1\le i\le r$; and
\item[(2)]
$h_t$ is not an identity of $L$,
\end{itemize}
where $d=d(L)$ is defined in (\ref{e5}).
\end{lemma}

{\em Proof.} Let $B_1,\ldots, B_m$ be as in Lemma \ref{l4}. Then by the definition
(up to reindexing of $G_1,\ldots,G_m$) there exist $b_1\in B_1,\ldots, b_k\in B_k$,
$a_1,\ldots,a_p\in L$ such that for some multilinear monomial $w(z_1,\ldots, z_{k+p})$
the value
$$
w(b_1,\ldots, b_k,a_1,\ldots,a_p)
$$
is non-zero and $\dim G_1+\cdots+\dim G_k=d$ in $G=L/N$.

Recall that for the adjoint representation of $G_i$ there exists a central polynomial
(see \cite[Theorem 12.1]{R}) i.e. an associative multilinear polynomial $g_i$ which 
assumes only scalar values on $ad~x_\alpha, x_\alpha \in G_i$. Moreover, $g_i$ is not 
an identity of the adjoint representation of $G_i$ and it depends on $q$ disjoint 
alternating sets of variables of order $d_i=\dim G_i$. That is,
\begin{equation}\label{e6}
g_i=g_i(x_{1,d_i}^1,\ldots,x_{d_i,d_i}^1,\ldots,x_{1,d_i}^q,\ldots,x_{d_i,d_i}^q)
\end{equation}
is skew-symmetric on each $\{x_{1,d_i}^j,\ldots,x_{d_i,d_i}^j\}$ and for some
$a_{1,d_i}^1,\ldots,a_{d_i,d_i}^q\in G_i$ the equality
$$
g_i(ad~a_{1,d_i}^1,\ldots,ad~a_{d_i,d_i}^q)(c_i)=c_i
$$
holds for any $c_i\in G_i$. Cleary, $a_{1,d_i}^j,\ldots, a_{d_i,d_i}^j$ are
linearly independent for any fixed $1\le j \le q$ since $g$ is alternating on
$\{x_{1,d_i}^j,\ldots, x_{d_i,d_i}^j\}$. Hence there exists an evaluation of
(\ref{e6}) in $L$ with all $\widetilde a_{\beta\gamma}^\alpha, \widetilde c_i$ 
in $B_i$ such that
\begin{equation}\label{e7}
g_i(ad~\widetilde a_{1,d_i}^1,\ldots,ad~\widetilde a_{d_i,d_i}^q)
(\widetilde c_i)\equiv \widetilde c_i({\rm mod}~N).
\end{equation}
On the other hand, if at least one of $\widetilde a_{\beta\gamma}^\alpha$ 
lies in $B_j$, $j\ne i$, or in $N$ then
\begin{equation}\label{e8}
g_i(ad~\widetilde a_{1,d_i}^1,\ldots,ad~\widetilde a_{d_i,d_i}^q)
(\widetilde c_i)\equiv 0~({\rm mod}~N).
\end{equation}
Since we can apply $g_i$ several times, the integer $q$ can be taken to be the same 
for all $i=1,\ldots,k$. Moreover, it follows from (\ref{e7}), (\ref{e8}) that for any 
$t=1,2,\ldots$ there exists a multilinear Lie polynomial
$$
f_i^t=f_i^t(x_{1,d_i}^1,\ldots,x_{d_i,d_i}^1,\ldots,x_{1,d_i}^{tq},
\ldots,x_{d_i,d_i}^{tq},y_i)
$$
alternating on each set $x_{1,d_i}^j,\ldots,x_{d_i,d_i}^j$, $1\le j\le tq$, such that
$$
f_i^t(\widetilde a_{1,d_i}^1,\ldots,\widetilde a_{d_i,d_i}^{tq},
\widetilde c_i)\equiv \widetilde c_i ({\rm mod}~N)
$$
for some $\widetilde a_{1,d_i}^1,\ldots,\widetilde a_{d_i,d_i}^{tq}\in B_i$ and 
for any $\widetilde c_i\in B_i$.

Recall that the monomial $w=w(z_1,\ldots,z_{k+q})$ has a non-zero evaluation
$$
\bar w=w(b_1,\ldots,b_k,a_1,\ldots,a_p)
$$
in $L$ with $b_1\in B_1,\ldots,b_k\in B_k$. Replacing $z_i$ by $f_i^t$ in $w$ for 
all $i=1,\ldots, k$ and alternating the result we obtain a polynomial
$$
h_t=Alt~w(f_1^t(x_{1,d_1}^1,\ldots,x_{d_1,d_1}^{tq},y_1),\ldots, 
f_k^t(x_{1,d_k}^1,\ldots,x_{d_k,d_k}^{tq},y_k), z_{k+1},\ldots, z_p),
$$
where $Alt=Alt_1\cdots Alt_{tq}$ and $Alt_j$ denotes the total alternation on variables
$$
x_{1,d_1}^j,\ldots,x_{d_1,d_1}^{j},\ldots,x_{1,d_k}^j,\ldots,x_{d_k,d_k}^{j}.
$$
Now if $\bar w=w(b_1,\ldots,b_k,a_1,\ldots,a_p)\in N^i\setminus N^{i+1}$ for some integer 
$i\ge 0$ in $L$ then according to (\ref{e7}), (\ref{e8}) we get 
$$
\rho(h_t)\equiv \left(d_1!\cdots d_k!\right)^{tq}\bar w ({\rm mod}~N^{i+1}),
$$
where $\rho:X\to L$ is an evaluation, $\rho(x^\alpha_{\beta\gamma})=a^\alpha_{\beta\gamma}$,
$\rho(y_j)=b_j$, $\rho(z_{k+j})=a_j$. In particular, $h_t$ is not an identity of $L$.
Renaming variables $x^\alpha_{\beta\gamma}, y_1,\ldots, y_k,z_{k+1},\ldots,z_{k+p}$
we obtain the required polynomial $h_t$ with $s=k+p$.

In order to get similar multialternating polynomial $h_t$ for $k+p+1$ we replace the 
initial polynomial $w=w(z_1,\ldots,z_{k+p})$ by $w'=w'(z_1,\ldots,z_{k+p+1}) =
w(z_1z_{k+p+1},z_2,\ldots,z_{k+p})$. Since $G_1$ is simple we have $G_1^2=G_1$. Hence
there exists an element $a_{p+1}\in B_1$ such that $w'(b_1,\ldots,b_k,a_1,\ldots,
a_{p+1})=w(b_1a_{p+1},b_2,\ldots,b_k, \\ a_1,\ldots,a_p)\ne 0$. Continuing this process we
obtain similar $h_t$ for all integers $k+p+2,k+p+3,\ldots~$.

\hfill $\Box$

Using multialternating polynomials constructed  in the previous lemma we get the following
lower bound for codimensions.

\begin{lemma}\label{l6}
Let $L,q$ and $s$ be as in Lemma \ref{l5}. Then there exists a constant $C>0$ such that
$$
c_n(L)\ge \frac{1}{Cn^{2d}}\cdot d^n
$$
for all $tq+s\le n\le tq+s+q-1$ and for all $t=1,2,\ldots~$.
\end{lemma}

{\em Proof}. Given $t$ and $s\le s'\le s+q-1$, consider the polynomial $h_t$ constructed
in Lemma \ref{l5}. Then $n=\deg h_t=tqd+s'$ and $h_t$ depends on $tq$ alternating sets of 
indeterminates of order $d$. Denote by $M$ the $FS_n$-submodule of $P_n(L)$ generated by
$h_t$. Let $n_0=tqd$ and let the subgroup $S_{n_0}\subseteq S_n$ act on $tqd$ alternating 
indeterminates $x_1^1,\ldots,x_d^{tq}$. Then $M_0=FS_{n_0}h_t$ is a nonzero subspace of 
$M$. Obviously,
\begin{equation}\label{e9}
c_n(L)\ge\dim M\ge \dim M_0.
\end{equation}

Consider character of $M_0$ and its decomposition onto irreducible components
\begin{equation}\label{e10}
\chi(M_0)=\sum_{\lambda\vdash n_0} m_\lambda\chi_\lambda.
\end{equation}
By Lemma \ref{l1} algebra $L$ satisfies Capelli identity of rank $d_0\ge \dim L/N\ge d$.
Hence $m_\lambda=0$ in (\ref{e10}) as soon as the height $ht(\lambda)$ of $\lambda$, i.e.
the number of rows in Young diagram $D_\lambda$, is bigger than $d_0$.

Now we prove that for any multilinear polynomial $f=f(x_1,\ldots,x_n)$ and for any partition 
$\lambda\vdash n_0$ with $\lambda_{d+1}\ge p$ where $N^p=0$ the polynomial $e_{T_\lambda}f$
is an identity of $L$.

Since $e_{T_\lambda}=R(T_\lambda)C(T_\lambda)$, it is sufficient to show that
$h=h(x_,\ldots, x_n)= C(T_\lambda) f$ is an identity of $L$. Note that the set 
$\{x_1,\ldots,x_{n_0}\}$ is a disjoint union
$$
\{x_1,\ldots,x_{n_0}\}=X_0\cup X_1\cup\ldots\cup X_p
$$
where $|X_1|,\ldots, |X_p|\ge d+1$ and $h$ is alternating on each 
$X_1,\ldots,X_p$, i.e. $h\in Q_{d+1,p}$. As it was shown in the proof of 
Lemma \ref{l4}, $h$ is an identity of $L$.

It follows that $m_\lambda\ne 0$ in (\ref{e10}) for $\lambda\vdash n_0$ only if
$\lambda_{d+1}<p$. In particular
\begin{equation}\label{e11}
n_0-(\lambda_1+\cdots+\lambda_d)\le(d-d_0)p.
\end{equation}

By the construction of essential idempotent $e_{T_\lambda}$ any polynomial 
$e_{T_\lambda}f(x_1,\ldots,x_{n_0})$ is symmetric on $\lambda_1$ variables corresponding
to the first row of $T_\lambda$. Since $h_t$ depends on $tq$ alternating sets of variables 
it follows that $e_{T_\lambda}h_t=0$ for any $\lambda\vdash n_0$ with $\lambda_1\ge tq+1$.

Denote $c_1=(d-d_0)p$. If $m_\lambda\ne 0$ in (\ref{e10}) for $\lambda\vdash n_0$,
$\lambda=(\lambda_1,\ldots,\lambda_k)$, then $k\le d_0$ and
\begin{equation}\label{e12}
\lambda_{d-1}\le\ldots\le \lambda_1\le tq.
\end{equation}
If $\lambda_d<tq-c_1$ then combining (\ref{e11}) and (\ref{e12}) we get
$\lambda_{d+1}+\cdots+ \lambda_k=n_0-(\lambda_1+\cdots+\lambda_d) \le c_1$ and
$$
n_0=(\lambda_1+\cdots+\lambda_{d-1})+\lambda_d+(\lambda_{d+1}+\cdots+\lambda_k) <
tq(d-1)+tq-c_1+c_1=tqd=n_0,
$$
a contradiction. Hence $\lambda_d\ge tq-c_1$ and Young diagram $D_\lambda$ contains
a rectangular diagram $D_\mu$ where
$$
\mu=(\underbrace{tq-c_1,\ldots,tq-c_1}_{d})
$$
is a partition of $n_1=d(tq-c_1)=n_0-c_1d = n-s'-c_1d \ge n-s-q-c_1d+1$ since 
$s'\le s+q-1$. From Hook formula for dimensions of irreducible $S_n$-representations
(see \cite[Proposition 2.2.8]{GZbook}) and from Stirling formula for factorials it 
easily follows that
$$
d_\mu=\deg \chi_\mu>\frac{d^{n_1}}{n_1^{2d}}
$$
for all $n$ sufficiently large. Since $\dim M_0\ge d_\lambda\ge d_\mu$ and $n_1\ge n-c_2$
for constant $c_2=s+q+c_1d-1$ we conclude from (\ref{e9}) that
$$
c_n(L)> \frac{d^{n}}{Cn^{2d}}
$$
where $C=d^{c_2}$ and we are done.

\hfill $\Box$

\section{Existence of PI-exponents}

It follows from Lemma \ref{l6} that
$$
\underline{\mbox{exp}}(L) = \liminf_{n\to \infty}\sqrt[n]{c_n(L)}\ge d.
$$
Combining this inequality with Lemma \ref{l4} we get the following

\begin{theorem}\label{t1}
Let $L$ be an almost nilpotent Lie algebra over a field $F$ of characteristic zero. If
$N$ is the maximal nilpotent ideal of $L$ and $L/N$ is semisimple then the PI-exponent
of $L$,
$$
\mbox{exp}(L)=\lim_{n\to \infty}\sqrt[n]{c_n(L)},
$$
exists and is a positive integer.
\end{theorem}
\medskip

Now consider the case of soluble almost nilpotent special Lie algebras.

\begin{theorem}\label{t2}
Let $L$ be an almost nilpotent soluble special Lie algebra over a field $F$ of characteristic zero. Then the PI-exponent
of $L$ exists and is a positive integer.
\end{theorem}
{\em Proof}. Let $L$ be a special soluble Lie algebra with a nilpotent ideal $N$ of a 
finite codimension. By Lemma \ref{l2} algebra $L$ satisfies Capelli identity of some rank.
Then the variety $V={\rm Var~}L$ generated by $L$ has a finite basis rank \cite{Z4}
that is $L$ has the same identities as some $k$-generated Lie algebra $L_k\in V$. Clearly,
$\underline{\mbox{exp}}(L)=\underline{\mbox{exp}}(L_k)$ and $\overline{\mbox{exp}}(L) =
\overline{\mbox{exp}}(L_k)$. Since $L$ is soluble it follows that $L_k$ is a finitely 
generated soluble Lie algebra from special variety $V$. By \cite[Proposition 6.3.2, 
Theorem 6.4.6]{Baht} we have $(L_k^2)^t=0$ for some $t \ge 1$. In this case
$\mbox{exp}(L)$ exists and is a non-negative integer \cite{PM}.

\hfill $\Box$


\begin{thebibliography}{99}

\bibitem{Baht}
Yu. A. Bahturin,
{\em Identical Relations in Lie algebras}, 
Utrecht, VNU Science Pres, 1987.

\bibitem{B-Dr} 
Yu. A. Bahturin and V. Drensky,
{\em Graded polynomial identities of matrices},  
Linear Algebra Appl.  {\bf 357}  (2002), 15--34.

\bibitem{GRZ1}
A. Giambruno, A. Regev and M. Zaicev,
{\em On the codimension growth of finite-dimensional Lie algebras}, 
J. Algebra {\bf 220} (1999), 466--474.

\bibitem{GRZ2}
A. Giambruno, A.  Regev and M. V. Zaicev, 
{\em Simple and semisimple Lie algebras and codimension growth}, 
Trans. Amer. Math. Soc.  {\bf 352}  (2000),  no. 4, 1935-1946.

\bibitem{GZ1}
A. Giambruno and M. Zaicev,
{\em On codimension growth of finitely generated associative 
algebras}, Adv. Math. {\bf 140} (1998), 145-155.

\bibitem{GZ2}
A. Giambruno and M. Zaicev,
{\em Exponential codimension growth of P.I.~algebras: an exact estimate}, 
Adv. Math. {\bf 142} (1999), 221-243.

\bibitem{GZbook} 
A. Giambruno and M. Zaicev, 
{\em Polynomial Identities and Asymptotic Methods},
Mathematical Surveys and Monographs Vol. {\bf 122}, 
American Mathematical Society, 
Providence, RI, 2005.

\bibitem{GZ-TAMS2010}
A. Giambruno and M. Zaicev,
{\em Codimension growth of special simple Jordan algebras}, 
Trans. Amer. Math. Soc. {\bf 362}  (2010), 3107--3123.

\bibitem{JK}
G. James and A. Kerber, 
{\em The Representation Theory of the Symmetric Group},
Encyclopedia of
Mathematics and its Applications, Vol.  {\bf 16}, Addison-Wesley, London, 1981.

\bibitem{M1}
S. P. Mishchenko, {\em On the problem of the Engel property}, (Russian),
Mat. Sb. (N.S.) {\bf 124 (166)} (1984), 56--67.

\bibitem{M2}
S. P. Mishchenko, 
{\em Growth of varieties of Lie algebras}, (Russian), Uspekhi Mat. Nauk.
{\bf 45} (1990), 25--45; 
English translation: Russian Math. Surveys {\bf 45} (1990), 27--52.

\bibitem{PM}
S. P. Mishchenko and P. M. Petrogradsky,
{\em Exponents of varieties of Lie algebras
with a nilpotent commutator subalgebra}, 
Comm. Alg. {\bf 27} (1999), 2223--2230.

\bibitem{P}
V. M. Petrogradsky, 
{\em Scale for codimension growth of Lie algebras},
Methods in Lie Theory (Levico Terme, 1997), 213--222, 
Lecture Notes in Pure and Appl. Math. {\bf 198} Marcel Dekker, New York, 1998.

\bibitem{R}
Yu. P. Razmyslov, {\em Identities of algebras and their representations},
Translations of Mathematical Monographs {\bf 138},
American Mathematical Society, Providence, RI, 1994.

\bibitem{Reg72}
A. Regev, {\em Existence of identities in $A \otimes B$}, 
Israel J. Math. {\bf 11} (1972), 131--152.

\bibitem{VMZ}
A. B. Verevkin, M. V. Zaicev and S. P. Mishchenko, 
{\em A sufficient condition for
coinciding of lower and upper exponents of a variety of linear algebras},  (Russian),
Vestnik Moscov. Univ. Ser. I Mat. Mekh. (to appear).

\bibitem{Z4}
M. V. Zaicev, {\em Standard identity in special Liealgebras},  (Russian),
Vestnik Moscow. Univ. Ser. I Mat. Mekh. No. 1 (1993), 56--59.

\bibitem{Z1}
M. Zaicev, {\em Identities of affine Kac-Moody algebras}, (Russian),
Vestnik Moscov. Univ. Ser. I Mat. Mekh. no. 2 (1996), 33--36; 
English translation:
Moscow Univ Math. Bull. {\bf 51} (1996), 29--31.

\bibitem{Z2}
M. Zaicev, {\em Varieties of affine Kac-Moody algebras}, (Russian),
Matem. Zametki {\bf 62} (1997), 95--102.

\bibitem{Z3}
M. Zaicev,  
{\em Integrality of exponents of growth of identities of finite-dimensional Lie algebras}, (Russian), Izv. Ross. Akad. Nauk Ser. Mat. {\bf 66} (2002), 23-48; 
English translation: Izv. Math. {\bf 66} (2002), 463-487.

\bibitem{ZM}
M. V. Zaicev and S. P. Mishchenko,
{\em An example of a variety of Lie
algebras with a fractional exponent}, 
J. Math. Sci. {\bf 93} (1999),
977--982.

\bibitem{ZM2}
M. V. Zaicev and S. P. Mishchenko, 
{\em On the polynomial growth of the colength of varieties
of Lie algebras}, (Russian), Algebra i Logika {\bf 38} (1999), 161--175; 
English translation:
Algebra and Logic {\bf 38} (1999), 84--92.


\end{thebibliography}
\end{document}